\documentclass[a4paper,10pt]{article}
\usepackage[affil-it]{authblk}
\usepackage[utf8]{inputenc}
\usepackage{amsmath,amssymb}
\usepackage{float}
\usepackage{array}
\usepackage{mathrsfs}
\usepackage{setspace}
\usepackage{subfig}
\usepackage{cancel}
\usepackage{graphicx}
\usepackage{txfonts}
\usepackage{anysize}
\usepackage{textcomp}
\usepackage{caption}
\usepackage{float}
\usepackage{epsfig}
\usepackage{color}
\title{An accurate and robust genuinely multidimensional Riemann solver for Euler equations based on TV flux splitting}
\author{S. Sangeeth  and  J. C. Mandal}
 \date{}
  \affil{Department of Aerospace Engineering,\\
Indian Institute of Technology Bombay,\\
Mumbai, 400076, India}
 
\begin{document}

\maketitle

\begin{abstract}A simple robust genuinenly multidimensional convective pressure split (CPS) , contact preserving, shock stable Riemann solver (GM-K-CUSP-X) for 
Euler equations of gasdynamics is developed. The convective and pressure components of the Euler system are seperated following the Toro-Vazquez
type PDE flux splitting [Toro et al, 2012].
Upwind discretization of these components are achieved using the framework of Mandal et al [Mandal et al, 2015].
The robustness of the scheme is studied 
on a few two dimensional test problems. The results demonstrate the efficacy of the scheme over the corresponding conventional two state version 
of the solver. Results from two classic strong shock test cases associated with the infamous Carbuncle phenomenon, indicate that the 
present solver is completely free of any such numerical instabilities albeit possessing contact resolution abilities.Such a finding 
emphasizes the preexisting notion about the positive effects that multidimensional flow modeling may have towards curing of shock instabilities.\newline
\end{abstract}

\section{Introduction}

Past few decades have witnessed commendable advancement in the computation of high speed flows. Birth of the upwind schemes for finite
volume methods is a 
milestone in this regard. Among these characteristics based schemes, those based on the solution of two state \textquoteleft Riemann problems\textquoteright at cell interface 
are the most popular class of schemes lately.
Godunov \cite{godunov_1959} proposed the first Riemann problem based algorithm (also called Riemann solvers) that used the exact solution of a Riemann problem to 
construct a first order accurate 
scheme. However, this exact solver was critized for being prohibitively expensive because of the iterative technique involved 
\cite{roe_1981}. 
To mitigate this shortcoming, a class of 
approximate Riemann solvers have been proposed and continue to be developed. Details about some famous approximate Riemann solvers can be found in 
\cite{roe_1981,harten_1983,toro_1994,liou_1993}. Such algorithms have brought forth an era of much cheaper yet
accurate computation of gasdynamic flows.

Since most of these approximate schemes are inherently one dimensional in their framework, they enjoy accurate and robust resolution of 
both linear and nonlinear wave fields arising in one dimensional problems. However, accurate resolution of flow features 
in practically relevant 
multidimensional problems where important flow features are oblique to the grid still pose a major challenge. This 
is mainly because the wave system in such problems possess 
infinitely many propagation direction as compared to the limited ones  
in a one-dimensional problem. Further vorticity 
wave enters the formulation and has to be dealt with appropriately \cite{roe_1986}.
Currently, the standard way of extending the one-dimensional schemes to multidimensions is through a dimensionally split operator
approach motivated in \cite{strang_1968} where local one dimensional Riemann problems are solved in a direction normal to each cell interface. 
Such a strategy has been criticised because of its unnatural selection of the interface normal as the only wave propagation direction even for problems
with infinitely many propagation directions.
This lacuna is attributed to cause for example, pressure disturbances across a grid oblique shear wave \cite{vanleer_1993}.

Another disturbing problem encountered by dimensionally split approximate Riemann solvers is the occurrence of various forms of numerical shock 
instabilities when simulating strong
normal shocks. A catalouge of many such failures can be found in \cite{quirk_1994}. It has been observed that only those schemes that have exact shear wave resolution
abilities are known to produce these instabilities. There is increasing evidence \cite{shen_2016,pandolfi_2001}
to consider that such failures occur due to lack of adequate dissipation across cell faces which are normal to the shock wave front. 
These problems solicit a strong need for a genuinely multidimensional formulation that resolves all characteristic fields while introducing 
appropriate dissipation along necessary directions.

Probably the earliest attempt at creating a multidimensional Riemann solver was by Raithby \cite{raithby_1976} who suggested 
discretizing the convective terms along flow dominant directions.
Davis \cite{davis_1984} and others \cite{levy_1993,nishikawa_2008,ren_2003} introduced the idea of 
rotated Riemann solvers that involves
identifying grid oblique shock orientations and solving Riemann problems across them. 


Roe and others \cite{roe_1986,rumsey_1991} suggested using extrapolated data from the left and right states across a Riemann solver 
to construct simple wave solutions for a multidimensional Riemann problem. This technique is not generally valid for arbitrary 
governing equations and demands substantial reformulation 
for three dimensional problems.


Collela \cite{collela_1990} is credited with proposing the CTU method that uses the characteristics of the system to incorporate diagonal cell contirbutions to interface
fluxes in a predictor corrector framework. 


Leveque \cite{leveque_1997} proposed a multidimensional scheme in which cross derivative terms that constitute the transverse contributions are 
included by interpreting the usual one dimensional Gudonov method in a fluctuation splitting framework. 


Ren et al \cite{ren_2006} constructed an operator split predictor corrector scheme based on CTU and Leveque's wave propagation method 
for both Euler and Navier Stokes systems. While the predictor step solves linearized Euler equations in characteristic variables, the corrector step adds viscous 
contributions.

Wendroff \cite{wendroff_1999} introduced a multistate Riemann solver that attempted the extension of one dimentional HLLE scheme \cite{harten_1983} into several dimensions. The 
scheme used an expensive nine point stencil and suffered unacceptable dissipation due to unnatural wavespeed selection \cite{balsara_2010}.

Another multistate Riemann solver is from the work of Brio et al \cite{brio_2001}. Multidimensional effects are incorporated by adding correction terms to the standard 
face normal fluxes at every computational cell interface. These correction terms are obtained by solving a three state Riemann problem using Roe's FDS at the corners 
of the respective cells. 

In a similar spirit, Balsara \cite{balsara_2010} presented a generic multidimensional HLLE solver (GM-HLLE) with simple closed form expressions 
that can be extended to any hyperbolic system. This method too relies on construction of multidimensional correction terms like \cite{brio_2001}
wherein these terms are obtained using a four state HLLE Riemann solver at interface corners.  

By building upon the wave model introduced in GM-HLLE, Balsara \cite{balsara_2012} further proposed a multidimensional HLLC solver for Euler and MHD systems. To deal 
with contact discontinuities in two dimensions, a set of twelve possible contact orientations on a given cell are included in the wave model.  
Although posessing closed form expressions in two dimensions, such a wave model would not be easily tractable in three dimensions.  

Improvements to the model was proposed in \cite{balsara_2014} by reformulating the scheme in terms of characteristic variables but this too remains 
complicated to implement on a computer code.

Mandal et al \cite{mandal_2015} proposed a multidimensional convective-pressure split scheme (GM-HLLCPS-Z) based on Zha-Bilgen \cite{zha_1993} way of
splitting Euler fluxes. The scheme consists of a wave speed averaged upwinding for the convective part and GM-HLLE type discretization for the pressure part.
This scheme is basically a multidimensional extension of HLL-CPS strategy \cite{mandal_2012} that uses a HLL type discretization 
for distinctly treating convective and pressure flux parts. The splitting of full Euler flux into convective and pressure parts can be 
achieved by adopting AUSM type or Zha-Bilgen type PDE splitting \cite{liou_1993,zha_1993}. Diffusion control is achieved by careful tuning of the 
dissipation vector of the pressure flux. This renders the scheme with stationary and moving contact preserving abilities in addition to 
capability of surviving the most stringent test problems. Toro has corroborated such a method by showing that exact contact preserving ability can be incorporated into convective-pressure
split framework by discretizing the pressure fluxes using a Riemann solver \cite{torovaz_2012}. Surprisingly, although contact preserving, HLL-CPS is found to 
evade most common forms of numerical shock instabilities, particularly the carbuncle phenomenon. 

Recently Xie et al \cite{xie2015_HLL} has proposed a contact capturing convective-pressure split scheme named K-CUSP-X. The scheme 
uses exactly similar discretization as HLL-CPS 
for both convective and pressure fluxes but differs only in the method of splitting the Euler fluxes into these components; instead of
using Zha-Bilgen or AUSM way of PDE
splitting, it adopts Toro-Vazquez method wherein the pressure terms embedded in the energy is also seperated \cite{torovaz_2012}. 
However unlike HLL-CPS, the present investigations reveal that this scheme is found to suffer from numerical shock instabilites.

In this paper, a new genuinely multidimensional scheme based on the conventional K-CUSP-X is proposed.  A multidimensional correction 
term similar to \cite{brio_2001} and \cite{balsara_2010} is
constructed by solving appropriate four state Riemann problem at the corners of each interface. The multidimensional terms are 
incorporated such that the final fluxes can be calculated at interfaces with the familiar ease of pre-existing dimensional split methods. 

Adhereing to the existing convention, this new scheme may be 
called GM-K-CUSP-X. It will be the purpose of this paper to demonstrate that such a construction is not only as accurate as GM-HLL-CPS-Z, 
but also cures the numerical
shock instability that plagued the corresponding two state conventional K-CUSP-X Riemann solver. The positive effect of genuinely multidimensional 
flow modeling on K-CUSP-X was first demonstrated in \cite{sangeeth_2016} wherein authors demonstrated that shock instability associated
with standing shock problem
\cite{dumbser_2004} was completely removed by extending the original two state solver into its genuinely multidimensional variant. Such a finding corroborates the pre existing opinion that 
multidimensional dissipation acts as a reliable cure for numerical shock instability problems and incentivizes the need for extending 
the accurate Riemann solvers in literature into their robust genuinely multidimensional counterparts.

This paper is organized as follows. In the next section, the governing equations and the type of PDE flux splitting used is detailed. In Section 3,
the second order version of the new formulation is described. In Section 4, results for some complex flow problems like double Mach reflection and
multidimensional Riemann problem are discussed. Further, two classic shock instability test problems, the odd-even decoupling problem and standing shock problem are
used to study instability behaviour of the newly developed scheme. Section 5 contains some concluding remarks.

\section{Preliminaries}
\subsection{Governing equation}
\label{gov_eq}
Consider the two-dimensional Euler equations in differential form
\begin{eqnarray}
\frac{\partial \mathbf{U}}{\partial t} + \frac{\partial \mathbf{F}}{\partial x} + \frac{\partial \mathbf{G}}{\partial y}=0
\label{euler}
\end{eqnarray}
where $\mathbf{U}=[\rho$ $\rho u$ $\rho v$ $\rho e]^T$ is the vector of conserved quantities. $\mathbf{F}$ and $\mathbf{G}$ are the flux vectors in the x and y directions respectively given by
\begin{eqnarray}
\mathbf{F}=\begin{bmatrix} \rho u \\ \rho u^2 + p \\ \rho uv \\ u(\rho e + p)\end{bmatrix} \hspace{10mm}
\mathbf{G}=\begin{bmatrix} \rho v \\ \rho uv \\ \rho v^2 + p \\ v(\rho e + p)\end{bmatrix}
\end{eqnarray}
In the present work, the above flux vectors are split into corresponding convective and pressure parts following 
the approach of Toro-Vazquez \cite{torovaz_2012} . Using ideal gas law, the split 
flux vectors can be written as
\begin{eqnarray}
\mathbf{F_1}=u\begin{bmatrix} \rho \\ \rho u \\ \rho v \\ \frac12 \rho (u^2+v^2) \end{bmatrix} \hspace{5mm}
\mathbf{F_2}=\begin{bmatrix} 0 \\ p \\ 0 \\ \frac{\gamma}{\gamma-1} pu \end{bmatrix} \hspace{10mm} 
\mathbf{G_1}=v\begin{bmatrix} \rho \\ \rho u \\ \rho v \\ \frac12 \rho (u^2+v^2) \end{bmatrix} \hspace{5mm}
\mathbf{G_2}=\begin{bmatrix} 0 \\ 0 \\ p \\ \frac{\gamma}{\gamma-1} pv \end{bmatrix}
\label{eqn:TVsplit}
\end{eqnarray}
where $\mathbf{F_1}$ and $\mathbf{G_1}$ are the convective fluxes while $\mathbf{F_2}$ and $\mathbf{G_2}$ are the pressure fluxes. $\gamma$ 
represents the ratio of specific heat capacities.
As suggested by Toro et al \cite{torovaz_2012} the convective fluxes $\mathbf{F_1}$ and $\mathbf{G_1}$ can be interpreted as 
simple advection of mass, momentum and kinetic energy along the x and y directions respectively,
while the pressure fluxes $\mathbf{F_2}$ and $\mathbf{G_2}$ are interpreted to be sonic impulses that spread in all directions 
with reference to these convecting particles. 
This type of PDE splitting primarily differs from those that already exist in the literature like \cite{zha_1993} and \cite{liou_1993}, 
in terms of the 
quantity concerning energy that is being advected. Such differences may have strong bearing on the robustness of these schemes. 

\subsection{Finite volume discretization}
\label{FVD}

	  \begin{figure}[H]
	  \centering
	  \includegraphics[scale=0.3]{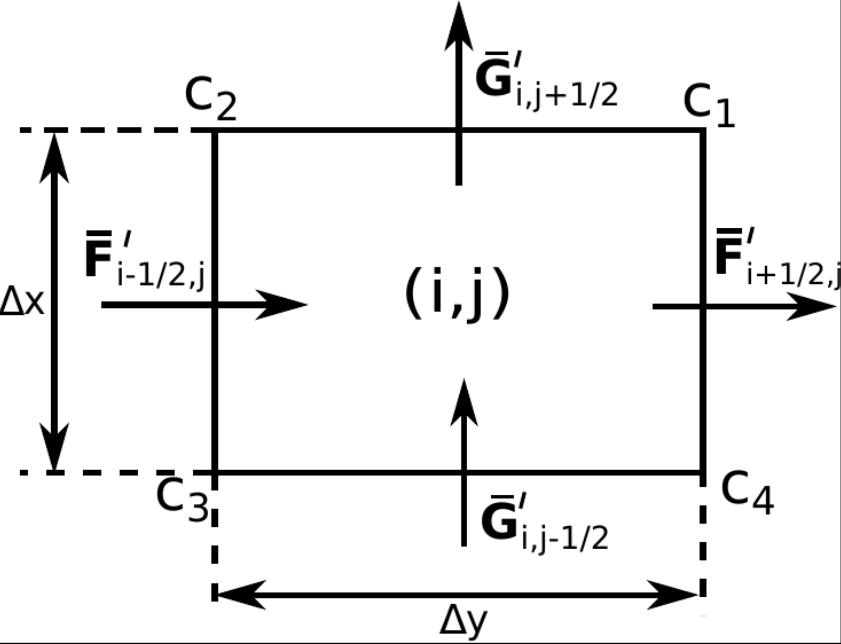}
	  \caption{Cell ($i,j$) of area $\Delta x \times \Delta y$ and corners $c_1$, $c_2$, $c_3$ and $c_4$ in a Cartesian grid.}
	  \label{fig:FVM}
	  \end{figure}
Consider the integral form of equation (\ref{euler})
\begin{eqnarray}
\frac{\partial}{\partial t}\int_\Omega \mathbf{U} d\Omega + \oint_{d\Omega}(\mathbf{F}+\mathbf{G})\cdot\hat{\mathbf{n}} ds=0
\end{eqnarray}
where $\hat{\mathbf{n}}=(n_x,n_y)$ is the unit vector along the face normal. Consider a Cartesian cell of area $\Delta x \times \Delta y$
as shown in Figure \ref{fig:FVM}. The finite volume discretization for the cell ($i,j$) can be written as
\begin{eqnarray}
\bar{\mathbf{U}}_{i,j}^{n+1}=\bar{\mathbf{U}}_{i,j}^n+\frac{\Delta t}{\Delta x}(\bar{\mathbf{F}}^{\prime}_{i-\frac12,j}-\bar{\mathbf{F}}^{\prime}_{i+\frac12,j})+\frac{\Delta t}{\Delta y}(\bar{\mathbf{G}}^{\prime}_{i,j-\frac12}-\bar{\mathbf{G}}^{\prime}_{i,j+\frac12})
\end{eqnarray}
where n is the time level, i and j are cell indices. The $\bar{(.)}$ quantities depict the respective averaged quantities. The total fluxes at the cell interfaces ($i+\frac12,j$) and ($i-\frac12,j$) are denoted as $\bar{\mathbf{F}}^{\prime}_{i+\frac12,j}$ and 
$\bar{\mathbf{F}}^{\prime}_{i-\frac12,j}$ respectively. Similarly $\bar{\mathbf{G}}^{\prime}_{i,j+\frac12}$ and 
$\bar{\mathbf{G}}^{\prime}_{i,j-\frac12}$ are the y directional total fluxes in corresponding y interfaces. 

\section{Formulation}

\begin{figure}[h]
	  \centering
	  \includegraphics[scale=0.3]{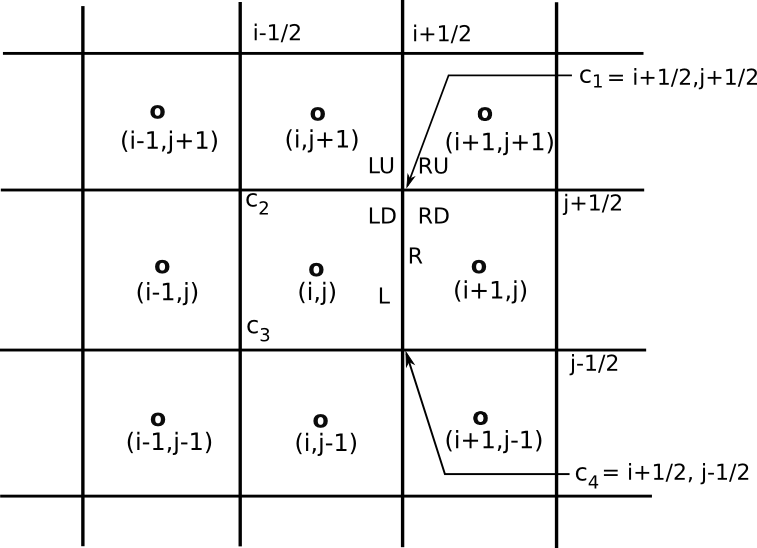}                
	  \caption{Figure showing a typical cell ($i,j$) in a Cartesian computational grid.}
	  \label{fig:gmgrid}
	  \end{figure}

\begin{figure}[h!]
	  \centering
	  \includegraphics[scale=0.3]{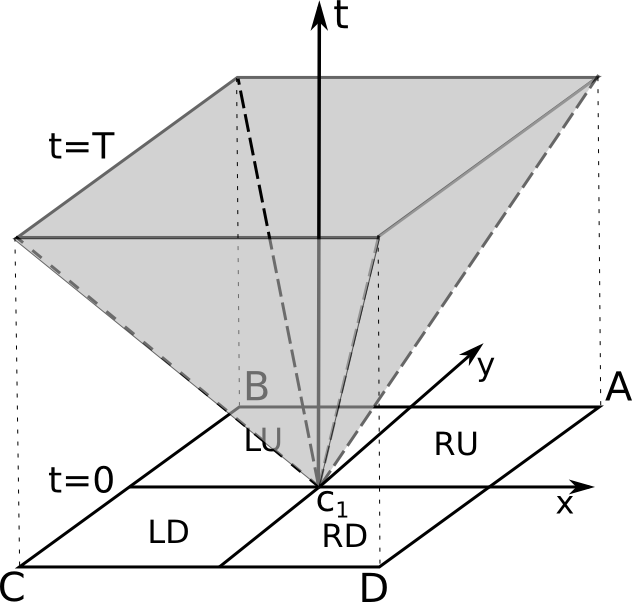}                
	  \caption{Figure showing the evolution of four waves making up the simple wave model used in this work.}
	  \label{fig:wavemodel}
	  \end{figure}

In the proposed scheme the interface flux will comprise of both a two state conventional Riemann flux and a  
multidimensional flux. These multidimensional
fluxes are sought from the solution of four state Riemann problem that occurs at each corner of a Cartesian cell. To see this more clearly, consider
a typical cell ($i,j$) as shown in Figure \ref{fig:gmgrid}. The four constant states that come together at corner $c_1$ at $t=0$ and forms 
a two-dimensional Riemann problem are marked as LU (Left-Upper), LD (Left-Down), RU (Right-Upper) and RD (Right-Down). Although realistically, 
the waves pertaining to this multistate Riemann problem will travel in infinitely many directions and swap a curvilinear area in 
space at any later time $t=T$, 
for sake of simplicity, the simple wave model consisting of only four waves as introduced in \cite{balsara_2010} is used in this work. 
Accordingly the wave propagation at the corners
is assumed to span, at any time T, a rectangular region. Figure \ref{fig:wavemodel} represents a three dimensional view of these 
waves as they evolve in time where the shaded region depicts the domain of influence of this multidimensional Riemann problem.

\begin{figure}[h]
	  \centering
	  \includegraphics[scale=0.3]{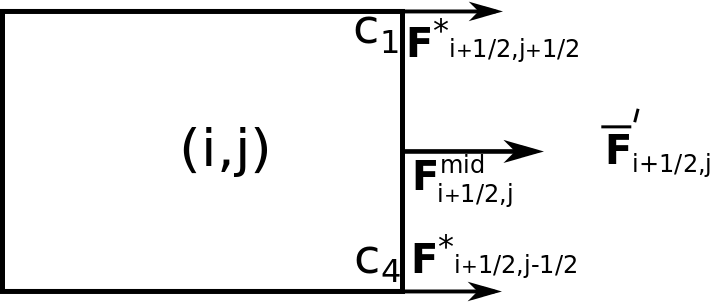}                
	  \caption{Figure showing the various components of an interface flux in the present genuinely multidimensional scheme.}
	  \label{fig:fluxDCU}
	  \end{figure}

In principle, the Equation \ref{euler} can be integrated along the space-time volume of Figure \ref{fig:wavemodel} appropriately 
to obtain the time averaged multidimensional
fluxes $\mathbf{F}^*$ in x-direction and $\mathbf{G}^*$ in y-direction. For a typical cell interface ($i+\frac12,j$) the total 
flux will be an ensemble of
the two multidimensional
fluxes from corners $c_1$ and $c_4$ denoted as $\mathbf{F}^*_{i+\frac12,j+\frac12}$, $\mathbf{F}^*_{i+\frac12,j-\frac12}$ and the 
single mid-point flux  $\mathbf{F}_{i+\frac12,j}^{mid}$. This is shown in Figure \ref{fig:fluxDCU}.
A conservative ensemble of the corner and mid point fluxes are achieved by using Simpson's rule of intergration
\cite{balsara_2010} along the interface as given by,

\begin{eqnarray}
\bar{\mathbf{F}}^{\prime}_{i+\frac12,j}=\frac16\mathbf{F}^*_{i+\frac12,j+\frac12}+\frac46\mathbf{F}_{i+\frac12,j}^{mid}+\frac16\mathbf{F}^*_{i+\frac12,j-\frac12}
\end{eqnarray}

\subsection{Evaluation of flux at the midpoint of the cell interface }
\label{midpoint_flux}

This section deals with determination of the two state Riemann flux at the mid point of a typical interface $(i+\frac12,j)$ denoted as 
${\mathbf{F}_{i+\frac12,j}^{mid}}$.  Following the Convective Pressure Split 
(CPS) philosophy the total flux at this interface is first split into convective and pressure parts. This work uses
the Toro-Vazquez type flux splitting as mentioned in Equation (\ref{eqn:TVsplit}),

\begin{eqnarray}
\mathbf{F}_{i+\frac12}^{mid} = {\mathbf{F}_1}_{i+\frac12,j}^{mid} + {\mathbf{F}_2}_{i+\frac12,j}^{mid}
\end{eqnarray}

These convective and pressure parts are discretized independently following the original HLL-CPS strategy \cite{mandal_2012}. It may be noted that the 
interface $(i+\frac12,j)$ admits only the x-directional Riemann flux ${\mathbf{F}_{i+\frac12,j}^{mid}}$ while the 
y-directional flux, ${\mathbf{G}_{i+\frac12,j}^{mid}}$ is zero on it. 


\subsubsection{Evaluation of convective flux at the midpoint of the cell interface (${\mathbf{F}_1}_{i+\frac12,j}^{mid}$)}

The upwind discretization of the convective flux denoted by ${\mathbf{F}_1}_{i+\frac12,j}^{mid}$ at the mid point of the cell interface ($i+\frac12,j$), 
following the strategy of HLL-CPS method \cite{mandal_2012}, is given by,
\begin{eqnarray}
{\mathbf{F}_1}_{i+\frac12,j}^{mid}=M_{k}\begin{bmatrix} \rho \\ \rho u \\ \rho v \\ \frac12 \rho (u^2+v^2)\end{bmatrix}_k a_k
\end{eqnarray}
\begin{eqnarray}
k=\begin{cases} L \quad \mbox{ if } \quad \bar{u}\geq0 \\ R \quad \mbox{ if } \quad \bar{u}<0\end{cases}
\end{eqnarray}
\begin{eqnarray}
M_{k}=\begin{cases} \displaystyle\frac{\bar{u}}{\bar{u}-S_L^c} \quad \mbox{ if } \quad \bar{u}\geq0 \\  \displaystyle\frac{\bar{u}}{\bar{u}-S_R^c} \quad \mbox{ if } \quad \bar{u}<0\end{cases}
\end{eqnarray}
\begin{eqnarray}
a_k=\begin{cases} u_{L}-S_L^c \quad \mbox{ if } \quad \bar{u}\geq0 \\ u_{R}-S_R^c \quad \mbox{ if } \quad \bar{u}<0\end{cases}
\end{eqnarray}
where, the average local x-directional velocity at the interface is taken as\, $\bar{u}=\displaystyle\frac{u_{L}+u_{R}}2$. Depending 
on the sign of the average local x-directional velocity, left (L) or right (R) states are selected for upwinding. 
$S_L^c$ and $S_R^c$ are carefully selected wave speeds which are discussed in section \ref{speedcenter}.

\subsubsection{Evaluation of pressure flux at the midpoint of the cell interface (${\mathbf{F}_2}_{i+\frac12,j}^{mid}$)}

The pressure flux at the midpoint of the cell interface is obtained by applying a HLL type discretization of the pressure flux vector 
\cite{mandal_2012}.
\begin{eqnarray}
\mathbf{F_2}_{i+\frac12,j}^{mid}=\frac{S_R^c}{S_R^c-S_L^c}{\mathbf{F}_2}_L-\frac{S_L^c}{S_R^c-S_L^c}{\mathbf{F}_2}_R+\frac{S_R^cS_L^c}{S_R^c-S_L^c}(\mathbf{U}_{R}-\mathbf{U}_{L})
\label{eqn:hll_p_flux}
\end{eqnarray}
The above equation can be rewritten as
\begin{eqnarray}
{\mathbf{F}_2}_{i+\frac12,j}^{mid}=\frac12(\mathbf{F}_{2L}+\mathbf{F}_{2R}) + \delta \mathbf{U}_2
\end{eqnarray}
where $\mathbf{F}_{2L}$ and $\mathbf{F}_{2R}$ are the left and right side x-directional pressure fluxes normal to the interface $(i+\frac12,j)$
and $\delta \mathbf{U}_2$ is the numerical diffusion given by
\begin{eqnarray}
\delta \mathbf{U}_2=\frac{S_R^c+S_L^c}{2(S_R^c-S_L^c)}({\mathbf{F}_2}_L-{\mathbf{F}_2}_R)-\frac{S_L^cS_R^c}{S_R^c-S_L^c}
\begin{bmatrix}\rho_{L}-\rho_{R}\\
               (\rho u)_{L}-(\rho u)_{R}\\
	       (\rho v)_{L}-(\rho v)_{R}\\
	       (\rho e)_{L}-(\rho e)_{R}
\end{bmatrix}
\label{eqn:pressureflux_hll_dissipationform}
\end{eqnarray}
It is clear from Equation \ref{eqn:pressureflux_hll_dissipationform} that the second term will give rise to numerical diffusion across a contact.
Thus in order to capture the contact discontinuity accurately, 
density terms in $\delta \mathbf{U}_2$ are replaced by the pressure terms by using isentropic assumption\, 
$\bar{a}_c^2=\displaystyle\frac{\delta p}{\delta \rho}$\, as described in the reference \cite{mandal_2012}.
\begin{eqnarray}
\delta \mathbf{U}_2=\frac{S_R^c+S_L^c}{2(S_R^c-S_L^c)}(\mathbf{F}_{2L}-\mathbf{F}_{2R})-\frac{S_R^cS_L^c}{\bar{a}_c^2(S_R^c-S_L^c)}
           \begin{bmatrix} p_L-p_R\\ (pu)_L-(pu)_R \\ (pv)_L-(pv)_R \\ \frac{\bar{a}_c^2}2 (p_L-p_R) + \frac12[(pq^2)_L-(pq^2)_R]\end{bmatrix}
\end{eqnarray} 
where $\bar{a}_c$ is the average speed of sound at the interface given by $\bar{a}_c=\displaystyle\frac{a_L+a_R}2$ and $q^2=u^2+v^2$ is twice the 
local kinetic energy per unit mass.

\subsubsection {Selection of wave speeds for the 1D Riemann problem at the interface ($\mathbf{i+\frac12,j}$) }\label{speedcenter}
The wave speeds are selected according to conventional HLL-CPS method \cite{mandal_2012}
\begin{align}
\nonumber
 S_L &= min (0, u_L-a_L,u^*-c^*)\\
 S_R &= max (0, u_r+a_r,u^*+c^*)
\end{align}
where $u^*$ and $c^*$ are given by \cite{mandal_2012},
\begin{align}
 \nonumber
 u^* &= \frac{u_L+u_R}{2}+\frac{a_L-a_R}{\gamma-1}\\
 c^* &= \frac{a_L+a_R}{2}+ \frac{\gamma-1}{4}(u_L-u_R)
\end{align}
Supersonic conditions are taken care by including '0' in the above expressions. For a stationary flow the wave speeds are modified as
\begin{eqnarray} 
S_L^c=-\bar{a}_c \hspace{10mm} S_R^c=\bar{a}_c
\end{eqnarray}
It is to be noted that midpoint fluxes at other interfaces can be obtained in the similar manner.

%

\subsection{Evaluation of flux at the corner of the cell interface}
\begin{figure}[h]
	  \centering
	  \includegraphics[scale=0.4]{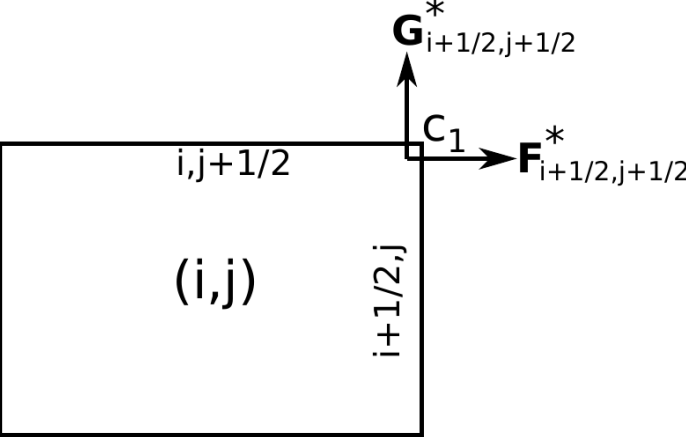}
	  \caption{Figure showing showing the multidimensional fluxes $\mathbf{F}^{*}_{i+\frac12,j+\frac12}$ and $\mathbf{G}^{*}_{i+\frac12,j+\frac12}$ 
	  at corner $c_1$ of $(i+\frac12,j)$ interface of a typical cell $(i,j)$ \cite{mandal_2015}.}
	  \label{fig:cornerflux}
	  \end{figure}

This section will detail how to evaluate the fluxes $\mathbf{F}^{*}_{i+\frac12,j+\frac12}$ and $\mathbf{G}^{*}_{i+\frac12,j+\frac12}$ that results 
from the interaction of four Riemann states at a representative corner $c_1$. These fluxes forms the 
multidimensional component of the interface fluxes and are shown in Figure \ref{fig:cornerflux}. Once again, resorting to the  
(CPS) philosophy, the total flux at this corner is split into convective ($\mathbf{{F_1}}^*_{i+\frac12,j+\frac12}$ and
$\mathbf{{G_1}}^*_{i+\frac12,j+\frac12}$) and 
pressure fluxes ($\mathbf{{F_2}}^*_{i+\frac12,j+\frac12}$ and $\mathbf{{G_2}}^*_{i+\frac12,j+\frac12}$) originating at the corner 
as per Equation \ref{eqn:TVsplit}. Analogous to the procedure for evaluation of the two state mid point flux at the interface, the split convective and the pressure parts at the 
corners will also be evaluated using different upwind strategies.

%

\subsubsection{Evaluation of convective flux at the corner $c_1$ of the cell interface  ($\mathbf{{F_1}}^*_{i+\frac12,j+\frac12}$)}

Following \cite{mandal_2015}, the x-directional convective flux is evaluated as, 

\begin{eqnarray}
\mathbf{{F_1}}^*_{i+\frac12,j+\frac12}=\bar{u}\frac{\left(S_U\begin{bmatrix} \rho \\ \rho u \\ \rho v \\ \frac12 \rho (u^2+v^2) \end{bmatrix}_{k_1}-S_D\begin{bmatrix} \rho \\ \rho u \\ \rho v \\ \frac12 \rho (u^2+v^2) \end{bmatrix}_{k_2}\right)}{S_U-S_D}
\label{eqn:convective_x}
\end{eqnarray}
where $\bar{u}$ is the wave speed averaged x-directional local fluid speed at the corner defined as, 
\begin{eqnarray}
\bar{u}=\frac{u_{LU}S_U-u_{LD}S_D+u_{RU}S_U-u_{RD}S_D}{2(S_U-S_D)}
\label{ubar_definition}
\end{eqnarray}
Based on the direction of the average x- 
directional flow, the upwind states $k_1,k_2$ are chosen as, \newline
If $\bar{u}>0$,\ $k_1=LU$ and $k_2=LD$ (i.e upwinding is done from the left states).\\
If $\bar{u}<0$,\ $k_1=RU$ and $k_2=RD$ (i.e upwinding is done from the right states).\\
In similar spirit, the y directional convective flux is evaluated as, 
\begin{eqnarray}
\mathbf{{G_1}}^*_{i+\frac12,j+\frac12}=\bar{v}\frac{\left(S_R\begin{bmatrix} \rho \\ \rho u \\ \rho v \\ \frac12 \rho (u^2+v^2) \end{bmatrix}_{k_1}-S_L\begin{bmatrix} \rho \\ \rho u \\ \rho v \\ \frac12 \rho (u^2+v^2) \end{bmatrix}_{k_2}\right)}{S_R-S_L}
\label{eqn:convective_y}
\end{eqnarray}
where $\bar{v}$ is the the wave speed averaged y-directional local fluid speed at the corner defined as,
\begin{eqnarray}
\bar{v}=\frac{v_{RU}S_R-v_{LU}S_L+v_{RD}S_R-v_{LD}S_L}{2(S_R-S_L)}
\label{vbar_definition}
\end{eqnarray}
The states $k_1,k_2$ are chosen accordingly as, \\
If $\bar{v}>0$,\ $k_1=RD$ and $k_2=LD$ (i.e upwinding is done from the down states).\\
If $\bar{v}<0$,\ $k_1=RU$ and $k_2=LU$ (i.e upwinding is done from the up states).\\
\par Since the above strategy is developed for a subsonic case, slight modification $\bar{u}$ and $\bar{v}$ is done to 
extend the above formulation to supersonic cases: \\
1. If the flow is supersonic in positive x-direction, then at the corner upwinding is done from left states. Therefore for the evaluation of x-directional convective flux $\bar{u}$ is taken as
\begin{eqnarray*}
\bar{u}=\frac{u_{LU}S_U-u_{LD}S_D}{S_U-S_D}
\end{eqnarray*}
2. If the flow is supersonic in negative x-direction the upwinding is done from right states. Therefore for the evaluation of x-directional  convective flux $\bar{u}$ is taken as
\begin{eqnarray*}
\bar{u}=\frac{u_{RU}S_U-u_{RD}S_D}{S_U-S_D}
\end{eqnarray*}
3. If the flow is supersonic in positive y-direction the upwinding is done from down states. Therefore for the evaluation of y-directional  convective flux $\bar{v}$ is taken as
\begin{eqnarray*}
\bar{v}=\frac{v_{RD}S_R-v_{LD}S_L}{S_R-S_L}
\end{eqnarray*}
4. If the flow is supersonic in negative y-direction the upwinding is done from upper states. Therefore for the evaluation of y-directional  convective flux $\bar{v}$ is taken as
\begin{eqnarray*}
\bar{v}=\frac{v_{RU}S_R-v_{LU}S_L}{S_R-S_L}
\end{eqnarray*}

\subsubsection{Evaluation of pressure flux at the corner $c_1$ of the cell interface ($\mathbf{{F_2}}^*_{i+\frac12,j+\frac12}$)}

Following \cite{mandal_2015}, the x-directional convective flux is evaluated as,
\begin{align}
\mathbf{{F_2}}^*_{i+\frac12,j+\frac12}&=\frac{\mathbf{{F_2}}_{LU}
S_RS_U+\mathbf{{F_2}}_{RD}S_LS_D-\mathbf{{F_2}}_{LD}S_RS_D-
\mathbf{{F_2}}_{RU}S_LS_U}{(S_R-S_L)(S_U-S_D)}  \nonumber \\ & -2\frac{S_RS_L}{(S_R-S_L)(S_U-S_D)}(\mathbf{{G_2}}_{RU}-\mathbf{{G_2}}_{LU}+\mathbf{{G_2}}_{LD}-
\mathbf{{G_2}}_{RD})  \nonumber \\ &+\frac{S_RS_L}{(S_R-S_L)(S_U-S_D)}(S_U(\mathbf{U}_{RU}-\mathbf{U}_{LU})-S_D(\mathbf{U}_{RD}-\mathbf{U}_{LD})) 
\label{eqn:balx}
\end{align}
Similarly the y-directional pressure flux is evaluated as, 
\begin{align}
\mathbf{G_2}^*_{i+\frac12,j+\frac12}&=\frac{\mathbf{{G_2}}_{RD}
S_RS_U+\mathbf{{G_2}}_{LU}S_LS_D-\mathbf{{G_2}}_{RU}S_RS_D-
\mathbf{{G_2}}_{LD}S_LS_U}{(S_R-S_L)(S_U-S_D)}  \nonumber \\ &-2\frac{S_US_D}{(S_R-S_L)(S_U-S_D)}(\mathbf{{F_2}}_{RU}-\mathbf{{F_2}}_{LU}+\mathbf{{F_2}}_{LD}-
\mathbf{{F_2}}_{RD})  \nonumber \\ &+\frac{S_US_D}{(S_R-S_L)(S_U-S_D)}(S_R(\mathbf{U}_{RU}-\mathbf{U}_{RD})-S_L(\mathbf{U}_{LU}-\mathbf{U}_{LD})) 
\label{eqn:baly}
\end{align}
The above equations can be rewritten as,
\begin{align}
\mathbf{F_2}^*_{i+\frac12,j+\frac12}&=\frac12(\mathbf{\overline{F}_2}
_L+\mathbf{\overline{F}_2}_R) \nonumber \\ &+ \delta\mathbf{U_{2x}}  -2\frac{S_RS_L}{(S_R-S_L)(S_U-S_D)}(\mathbf{{G_2}}_{RU}-\mathbf{{G_2}}_{LU}+\mathbf{{G_2}}_{LD}-
\mathbf{{G_2}}_{RD})
\label{hllbalx}
\end{align}

\begin{align}
\mathbf{G_2}^*_{i+\frac12,j+\frac12}&=\frac12(\mathbf{\overline{G}_2}
_D+\mathbf{\overline{G}_2}_U) \nonumber \\ &+ \delta\mathbf{U_{2y}} -2\frac{S_US_D}{(S_R-S_L)(S_U-S_D)}(\mathbf{{F_2}}_{RU}-\mathbf{{F_2}}_{LU}+\mathbf{{F_2}}_{LD}-
\mathbf{{F_2}}_{RD})
\label{hllbaly}
\end{align}
where 

\begin{eqnarray}
\mathbf{\overline{F}_2}_L=\frac{\mathbf{F}_{LU}S_U-\mathbf{F}_{LD}S_D}{S_U-S_D}
\end{eqnarray} 

\begin{eqnarray}
\mathbf{\overline{F}_2}_R=\frac{\mathbf{F}_{RU}S_U-\mathbf{F}_{RD}S_D}{S_U-S_D}
\end{eqnarray}

\begin{eqnarray}
\mathbf{\overline{G}_2}_D=\frac{\mathbf{G}_{RD}S_R-\mathbf{G}_{LD}S_L}{S_R-S_L}
\end{eqnarray}

\begin{eqnarray}
\mathbf{\overline{G}_2}_U=\frac{\mathbf{G}_{RU}S_R-\mathbf{G}_{LU}S_L}{S_R-S_L}
\end{eqnarray}
The  $\delta\mathbf{U_{2x}}$ and $\delta\mathbf{U_{2y}}$ terms in Equations (\ref{hllbalx},\ref{hllbaly}), are the 
numerical diffusion terms in x and y-directions respectively. To remove the numerical dissipation across a contact wave, the dissipation
terms in x and y-directions are remodeled using isentropic expression as, 

\begin{align}
\delta\mathbf{U_{2x}}&=\frac{S_R+S_L}{2(S_R-S_L)}(\mathbf{\overline{F}_2}
_L-\mathbf{\overline{F}_2}_R)\nonumber \\ &-\frac{S_RS_L}{(S_R-S_L)(S_U-S_D)(\bar{a}^2)}\begin{bmatrix}S_U(p_{LU}-p_{RU})-S_D(p_{LD}-p_{RD})\\
										  S_U((pu)_{LU}-(pu)_{RU})-S_D((pu)_{LD}-(pu)_{RD})\\
											S_U((pv)_{LU}-(pv)_{RU})-S_D((pv)_{LD}-(pv)_{RD})\\
											S_U(e^*_{LU}-e^*_{RU})-S_D(e^*_{LD}-e^*_{RD})
											\end{bmatrix}
\label{diffx}
\end{align}

\begin{align}
\delta\mathbf{U_{2y}}&=\frac{S_U+S_D}{2(S_U-S_D)}(\mathbf{\overline{G}_2}_D-\mathbf{\overline{G}_2}_U)\nonumber \\ &-\frac{S_US_D}{(S_R-S_L)(S_U-S_D)(\bar{a}^2)}\begin{bmatrix}S_R(p_{RD}-p_{RU})-S_L(p_{LD}-p_{LU})\\
											S_R((pu)_{RD}-(pu)_{RU})-S_L((pu)_{LD}-(pu)_{LU})\\
											S_R((pv)_{RD}-(pv)_{RU})-S_L((pv)_{LD}-(pv)_{LU})\\
											S_R(e^*_{RD}-e^*_{RU})-S_L(e^*_{LD}-e^*_{LU})
											\end{bmatrix}
\label{diffy}
\end{align}
where $e^*_k$ is given as
\begin{eqnarray*}
e^*_k=\frac{\bar{a}^2}{\gamma-1}p_{k}+\frac12p_k(u^2+v^2)_{k}
\end{eqnarray*}
and $\bar{a}$ as
\begin{eqnarray*}
\bar{a}=\frac{a_{LU}+a_{RU}+a_{LD}+a_{RD}}4
\end{eqnarray*}

The flux contributions due to the genuinely multidimensional Riemann problem at the other corner $c_4$ of the interface 
can be obtained in a similar manner. Most importantly, it must be noted that
while $\mathbf{F}^*_{i+\frac12,j \pm \frac12}$ 
contributes to the total interface flux at the interface $(i+\frac12,j)$, $\mathbf{G}^*_{i + \frac12,j + \frac12}$ and 
$\mathbf{G}^*_{i + \frac12,j - \frac12}$ contributes to the 
total interface flux at the interfaces $(i,j+\frac12)$ and $(i,j-\frac12)$ respectively.
%
\subsection{Selection of wave speeds for the multidimensional Riemann problem at the corners }

\begin{figure}[h]
	  \centering
	  \includegraphics[scale=0.4]{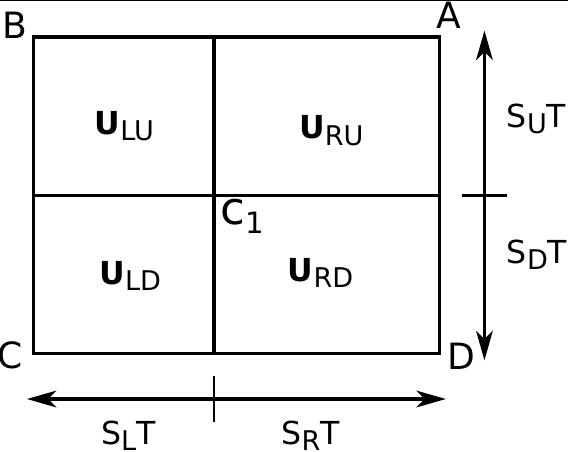}
	  \caption{Figure showing the rectangular area swept by the four waves that constitute the simple wave model used in this work at any time t=T. \cite{mandal_2015}}
	  \label{fig:wavemodelprojection}
	  \end{figure}

As previously mentioned, the present work adopts a simple wave model as proposed by \cite{balsara_2010} to represent the waves emerging from 
the four state Riemann problem at the corners of every interface.
A top view of the area swept by these four waves $S_L,S_R,S_U,S_D$ is shown in Figure \ref{fig:wavemodelprojection}. The rectangle ABCD in Figure \ref{fig:wavemodelprojection} 
depicts the domain of influence of the four state Riemann problem at corner $c_1$ at time T on x-y plane.
An estimate for these wavespeeds can be obtained as,
\begin{align}
\nonumber
S_R&=max(0,\lambda_x^N(\mathbf{U}_{RU}),\lambda_x^N(\mathbf{U}_{RD}),\bar{\lambda}_x^N(\mathbf{U}_{LU},\mathbf{U}_{RU}),\tilde{\lambda}_x^N(\mathbf{U}_{LD},\mathbf{U}_{RD}))\\
\nonumber
S_L&=min(0,\lambda_x^1(\mathbf{U}_{LU}),\lambda_x^1(\mathbf{U}_{LD}),\bar{\lambda}_x^1(\mathbf{U}_{LU},\mathbf{U}_{RU}),\tilde{\lambda}_x^1(\mathbf{U}_{LD},\mathbf{U}_{RD}))\\
\nonumber
S_U&=max(0,\lambda_y^N(\mathbf{U}_{RU}),\lambda_y^N(\mathbf{U}_{LU}),\bar{\lambda}_y^N(\mathbf{U}_{RD},\mathbf{U}_{RU}),\tilde{\lambda}_y^N(\mathbf{U}_{LD},\mathbf{U}_{LU}))\\
\nonumber
S_D&=min(0,\lambda_y^1(\mathbf{U}_{RD}),\lambda_y^1(\mathbf{U}_{LD}),\bar{\lambda}_y^1(\mathbf{U}_{RD},\mathbf{U}_{RU}),\tilde{\lambda}_y^1(\mathbf{U}_{LD},\mathbf{U}_{LU}))\\
\end{align}
where typically,\\
$\lambda_x^1(\mathbf{U}_{k})$ denote smallest x-directional wave speed in the state $\mathbf{U}_{k}$,\\
$\lambda_x^N(\mathbf{U}_{k})$ denote largest x-directional wave speed in the state $\mathbf{U}_{k}$,\\
$\tilde{\lambda}_x^1(\mathbf{U}_{k},\mathbf{U}_{l})$ denotes smallest x-directional wave speed from Roe averaged state between $\mathbf{U}_{k}$ and $\mathbf{U}_{l}$,\\
$\tilde{\lambda}_x^N(\mathbf{U}_{k},\mathbf{U}_{l})$ denotes largest x-directional wave speed from Roe averaged state between $\mathbf{U}_{k}$ and $\mathbf{U}_{l}$
such that $k,l \in \{LU,LD,RU,RD\}$ and $k\neq l$.\\
If $\bar{u}=0$ and $\bar{v}\neq0$, then x-directional wave speeds are modified as $S_L=-\bar{a}$ and $S_R=\bar{a}$.\\
If $\bar{u}\neq0$ and $\bar{v}=0$, then y-directional wave speeds are modified as $S_D=-\bar{a}$ and $S_U=\bar{a}$.\\
If $\bar{u}=0$ and $\bar{v}=0$, then wave speeds are modified as $S_L=-\bar{a}$, $S_R=\bar{a}$,  $S_D=-\bar{a}$ and $S_U=\bar{a}$.\\
It should be noted that zero has been added in the above expressions in order to ensure fully one sided flux in supersonic flow.

\section{Results}

\subsection{Isentropic vortex problem}
A second order accurate version of the present solver is developed using the SDWLS strategy \cite{mandal_2011} whose
details are omitted here for brevity.
Formal order of accuracy of the second order version of GM-K-CUSP-X is investigated using the isentropic vortex problem \cite{mandal_2015}.  
The problem involves an isentropic vortex centered initially at (0,0) and made to traverse the domain diagonally
under a periodic boundary condition. A Cartesian domain of size [-5,5] $\times$ [-5,5] is used. The initial conditions consists of a
unperturbed state given by $(\rho,p,u,v)$= $(1.0,1.0,1.0,1.0)$. The temperature is defined as $T=\frac{p}{\rho}$. A perturbation is
added to the flow as, 

\[
(\delta u, \delta v) =\displaystyle\frac{\epsilon}{2\pi}e^{0.5(1-r^2)}(-y,x)
\]
\[
\delta T=-\displaystyle\frac{\gamma-1}{8 \gamma \pi^2}\epsilon^2 e^{(1-r^2)}
\]
\[
\delta \rho=\delta T^{\displaystyle\frac{1}{\gamma-1}}
\]
\[
\delta p=\delta\rho^\gamma
\]
Here, $\epsilon$ which defines the strength of the vortex is taken as 5.0. $r$ denotes the Cartesian distance from vortex center. The 
accuracy of the scheme is measured in $L_1$ and $L_{\infty}$ norms of the density variable. The formal order of accuracy $O$ is obtained by
the formula, 
\begin{align}
 O &= \frac{log_{10}(\eta_2) - log_{10}(\eta_1)}{log_{10}(\Delta x_2) - log_{10}(\Delta x_1)}
\end{align}
where, $\eta_1$ and $\eta_2$ depict consecutive norms ($L_1$ or $L_{\infty}$) for progressively refined grids with dimensions
$\Delta x_1$ and $\Delta x_2$ respectively. The results are tabulated in table \ref{table:secondorderaccuracygmkcuspx}.

\begin{table}[!ht]
\centering 
\begin{tabular}{|c|c|c|c|c|} 
\hline
Mesh size & $L_1$ error & $L_1$ order & $L_{\infty}$ error & $L_{\infty}$ order\\ [0.5ex] 
\hline 
64$\times$64 & 0.007724 & - & 0.145543  & - \\ 
128$\times$128 & 0.001949 & 1.9866 & 0.044705  & 1.7029\\
256$\times$256 & 0.000510 & 1.9341 & 0.013459  & 1.7318 \\
512$\times$512 & 0.000110 & 2.2129 & 0.002416  & 2.4778 \\ [1ex] 
\hline 
\end{tabular}
\caption{Second order accuracy analysis of GM-K-CUSP-X using SDWLS reconstruction strategy.} 
\label{table:secondorderaccuracygmkcuspx} 
\end{table}
It is observed from the analysis that GM-K-CUSP-X scheme is able to achieve second order accuracy on sufficiently refined grids in
both $L_1$ and $L_{\infty}$ norms.

\subsection{Two dimensional Riemann problems}

Two dimensional Riemann problems provides an excellent test case to assess the qualitative improvement of genuinely multidimensional formulation
of GM-K-CUSP-X over the corresponding conventional two state K-CUSP-X scheme. Second order versions of both solvers are used for these test cases.
The first two dimensional Riemann problem investigated in the present work consists of a double Mach reflection and an oblique shock 
wave propagating at an angle to the grid.  
The initial conditions of the 2D Riemann problem are given by \cite{mandal_2015},

\begin{table}[!ht]
\centering 
\begin{tabular}{|c|c|c|c|c|} 
\hline
Zone & $\rho$ & $p$ & $u$ & $v$\\ [0.5ex] 
\hline 
$x>0$,\ $y>0$ & 1.5 & 1.5 & 0.0 & 0.0  \\ 
$x>0$,\ $y<0$ & 0.5323 & 0.3 & 0.0 & 1.206  \\
$x<0$,\ $y>0$ & 0.5323 & 0.3 & 1.206 & 0.0  \\
$x<0$,\ $y<0$ & 0.1379 & 0.029 & 1.206 & 1.206  \\
\hline 
\end{tabular}
\caption{Initial condition for Riemann problem 1.} 
\label{table:rp1_2000} 
\end{table}

\begin{figure}[H]
	  \centering
	  \subfloat[GM-K-CUSP-X]{\includegraphics[scale=0.5]{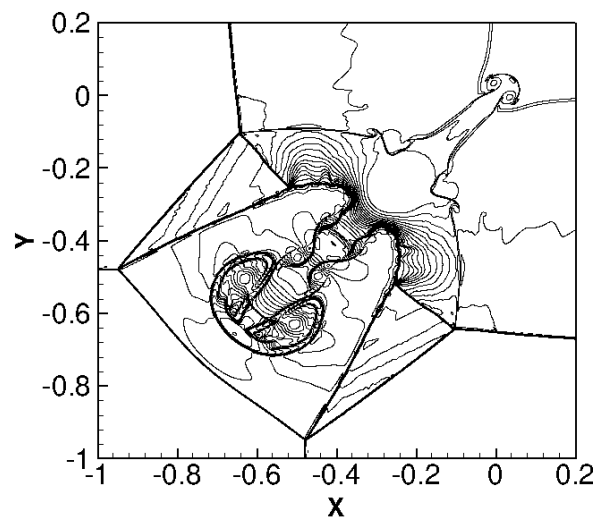}}\\
	  \subfloat[Original K-CUSP-X]{\includegraphics[scale=0.5]{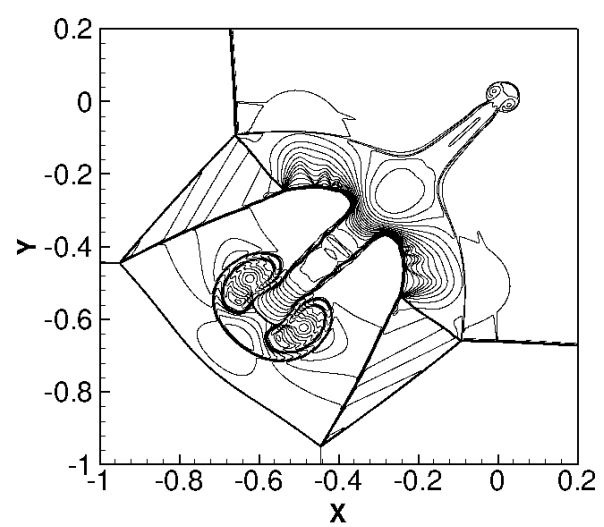} }
	   \caption{ Results for the Riemann problem 1 showing 30 iso density contours from 0.15 to 1.7 }
	   \label{rp1_2000}
	  \end{figure}

A Cartesian grid of size $2000\times2000$ spanning $[-1,1]\times[-1,1]$ is used.
The solution is evolved up to a time of 1.05 units with a CFL number of 0.95 as suggested in reference \cite{balsara_2012}. 
The density contours at the final time are shown in Figure \ref{rp1_2000}. Result for K-CUSP-X solver
under identical conditions is given for comparison. It is observed that the mushroom cap structure is well resolved by GM-K-CUSP-X as compared
to the original K-CUSP-X. Further GM-K-CUSP-X is able to resolve the Kelvin-Helmholtz roll up much better than K-CUSP-X. 
The second Riemann problem consist evolution of two weak shock waves and two contact waves and
can be set using conditions \cite{mandal_2015},

\begin{table}[!ht]
\centering 
\begin{tabular}{|c|c|c|c|c|} 
\hline
Zone & $\rho$ & $p$ & $u$ & $v$\\ [0.5ex] 
\hline 
$x>0$,\ $y>0$ & 0.5313 & 0.4 & 0.0 & 0.0  \\ 
$x>0$,\ $y<0$ & 1.0 & 1.0 & 0.0 & 0.7276  \\
$x<0$,\ $y>0$ & 1.0 & 1.0 & 0.7276 & 0.0  \\
$x<0$,\ $y<0$ & 0.8 & 1.0 & 0.0 & 0.0  \\
\hline 
\end{tabular}
\caption{Initial condition for Riemann problem 2.} 
\label{table:rp3_2000} 
\end{table}

\begin{figure}
	  \centering
	  \subfloat[GM-K-CUSP-X]{\includegraphics[scale=0.48]{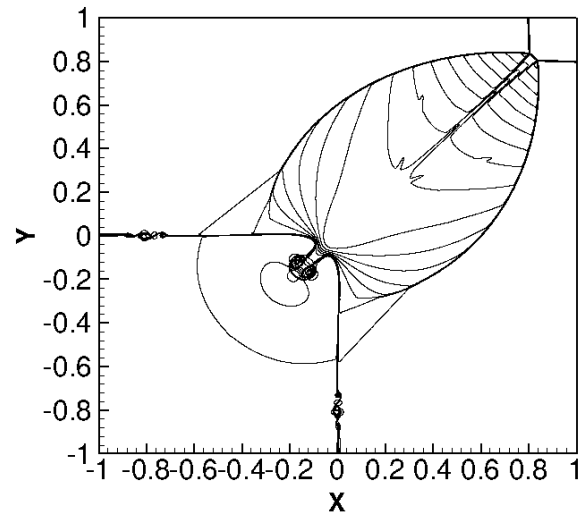}} \\
	  \subfloat[Original K-CUSP-X]{\includegraphics[scale=0.43]{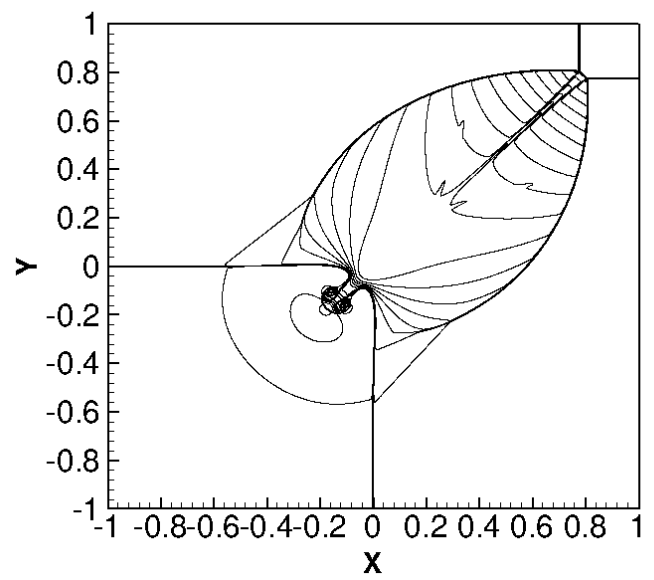} }
	  \caption{ Results for the Riemann problem 2 showing 30 iso density contours from 0.5 to 1.7 }
	  \label{rp3_2000}
	  \end{figure}

A 2000X2000 Cartesian grid spanning [-1,1]$\times$[-1,1] domain is used. CFL number of 0.95 is used and the solution is evolved for a
time of 0.5 units. The results obtained are represented using iso density contours. For comparison, results 
obtained under identical conditions by corresponding two state conventional K-CUSP-X is also provided in Figure \ref{rp3_2000}. 
Although both solvers are able to resolve the resultant contact waves and Mach reflections \cite{balsara_2010}, only GM-K-CUSP-X is able to 
resolve the Kelvin Helmholtz instability along the Mach stems.
	  
\subsection{Double Mach reflection problem}

This problems deals with a Mach 10 shock inclined at $60^o$ with the x-axis and propagating downstream of a rectangular duct of size 
[0,4]$\times$[0,1] and interacting with a reflective bottom boundary wall. Detailed initial and boundary conditions can be found in 
\cite{mandal_2015}. The problem is solved on a Cartesian mesh of size 1920$\times$480 with a CFL of 0.7 and evolved up to time of 0.2
units. The second order accurate results obtained for GM-K-CUSP-X are represented using iso density contours and are compared with that 
obtained using K-CUSP-X under same conditions. It is visible from the result that the present solver is able to resolve the complex shock structures
crisply and also the slipping contact line emerging from the triple point. 

	  \begin{figure}[H]
	  \centering
	  \subfloat[GM-K-CUSP-X]{\includegraphics[scale=0.42]{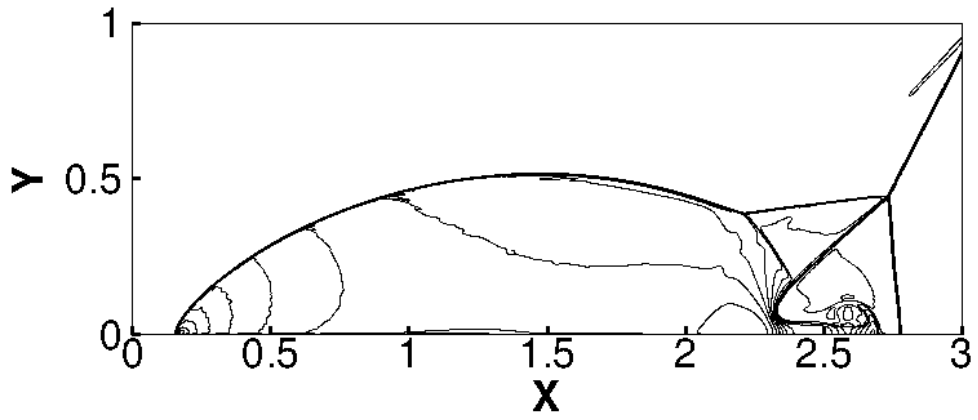}}\\
	  \subfloat[Original K-CUSP-X]{\includegraphics[scale=0.42]{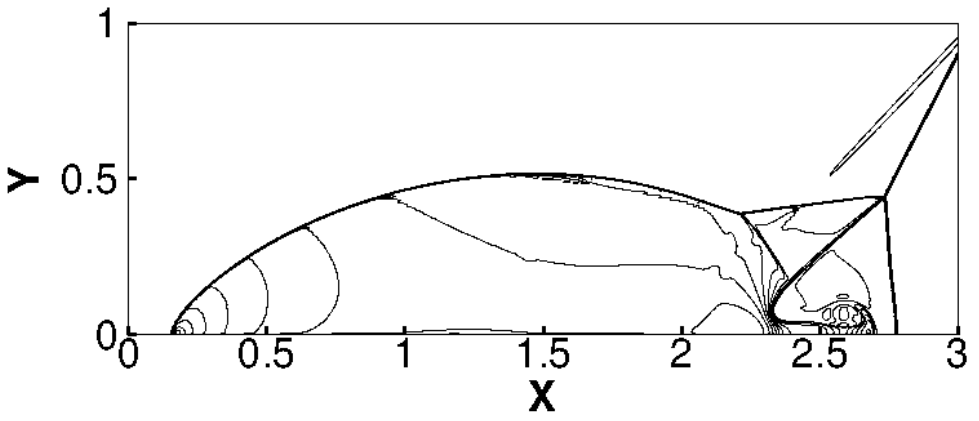} }
	  \caption{ Results for double Mach reflection problem showing 25 iso density contours from 1.77 to 22.44 }
	  \end{figure}

\section{Numerical shock instability tests}

One of the objectives of the present work is to demonstrate the effect 
of genuinely multidimensional formulation on the numerical shock instability 
characteristics of a Riemann solver. Two standard test cases namely, the odd-even decoupling problem \cite{quirk_1994} and the standing 
shock problem \cite{dumbser_2004} will be used to carry out the investigations. Since shock instability is most explicitly observed in the first
order simulations, first order version of solvers will be used for these tests. 

\subsection{Odd-Even decoupling problem}
\begin{figure}[H]
	  \centering
	  \setcounter{subfigure}{0}
	  \subfloat[K-CUSP-X]{\label{fig:kcusp_oddeven}\includegraphics[width=10.85cm,height=2.85cm]{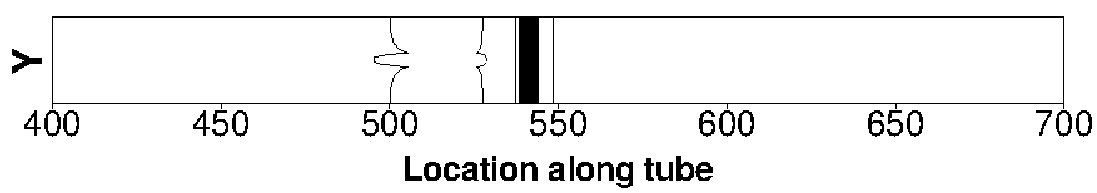}}\\
	  \subfloat[GM-K-CUSP-X]{\label{fig:gmkcusp_oddeven}\includegraphics[scale=0.4]{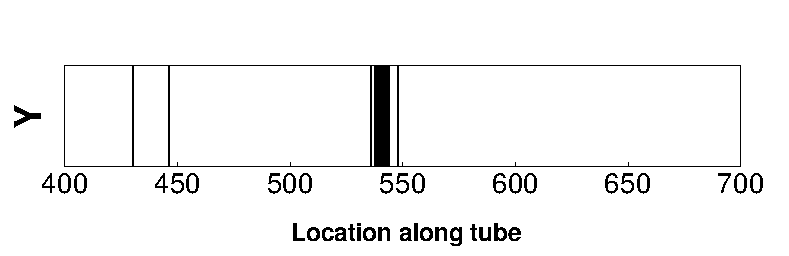} }
	 \caption{ Results for odd-even decoupling problem showing 40 iso density contours from 1.4 to 7.37 }
	  \end{figure}

Odd-even decoupling test consists of a M = 6 shock propagating down a computational duct whose centerline grid is slightly perturbed 
to induce random oscillations into the initial conditions \cite{quirk_1994}. It has been long argued that odd-even decoupling and the Carbuncle have the
same origin \cite{Chauvat_2005}. Much alike the Carbuncle phenomenon, the moving shock profile in the odd-even decoupling problem 
also deteriorates over time (producing the recognizable bulge at the centerline) and pollutes the after shock flow field. The original 
K-CUSP-X scheme was found to have slight after shock perturbations in this test. In comparison, the behavior 
of GM-K-CUSP-X scheme on this problem is shown in \ref{fig:gmkcusp_oddeven} after the solution has evolved for t = 140 units.
It is clearly observed that the genuinely multidimensional extension is able to preserve the shock profile without any instabilities.

\subsection{Standing shock instability}
\begin{figure}[H]
	  \centering
	  \setcounter{subfigure}{0}
	  \subfloat[K-CUSP-X]{\label{fig:kcusp_standingshock}\includegraphics[scale=0.3]{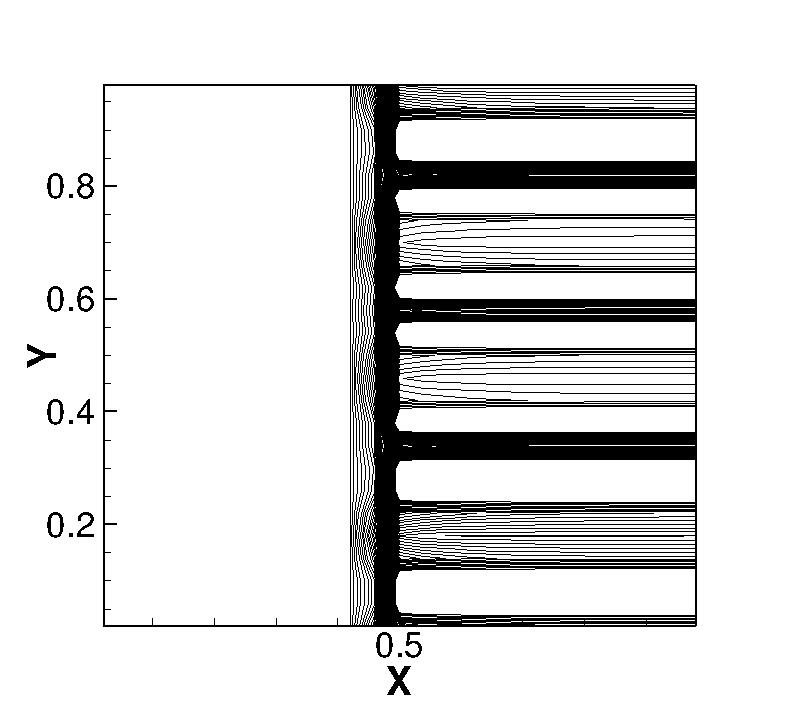}}\\
	  \subfloat[GM-K-CUSP-X]{\label{fig:gmkcusp_standingshock}\includegraphics[scale=0.3]{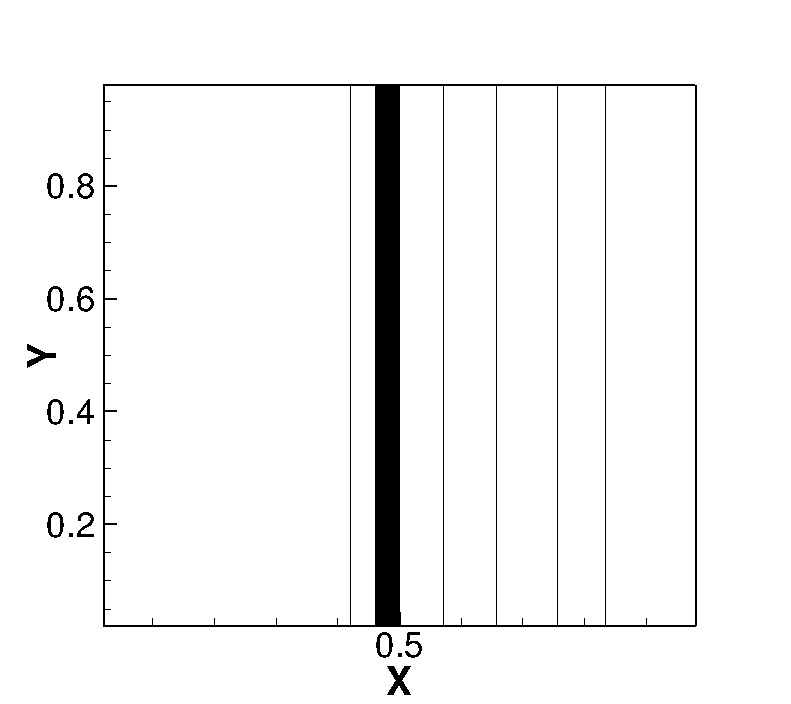} }
	 \caption{ Results for odd-even decoupling problem showing 50 iso density contours from 1.0 to 5.4 }
	  \end{figure}

To clearly differentiate the behavior of K-CUSP-X and GM-K-CUSP-X schemes, a standing shock instability test case is used. The details 
of the test are available in \cite{dumbser_2004}. As seen in the figure \ref{fig:kcusp_standingshock}, K-CUSP-X fails miserably in this
test case. However, from figure \ref{fig:gmkcusp_standingshock} it is very evident that GM-K-CUSP-X scheme is free of this instability.   

The reason for the stabilizing nature of GM-K-CUSP-X can be discerned from observing Equations (\ref{eqn:convective_x}), (\ref{eqn:convective_y}), 
(\ref{hllbalx}), (\ref{hllbaly}). Discretization of the convective terms using Equations (\ref{eqn:convective_x}), (\ref{eqn:convective_y}) 
employs a wave weighted averaging that strives to accurately model the underlying multidimensional wave evolution phenomenon. Specifically, it is
easy to note that Equations (\ref{ubar_definition}) and (\ref{ubar_definition}) defines wave averaged convection velocities in x and y directions
that depends on all the adjoining states at the corners. Such a formulation admits information from transverse cells which would contribute to 
dissipation that would help supressing the instabilities.
A similar observation can be found from Equations (\ref{hllbalx}), (\ref{hllbaly}) which describe the discretization of the pressure terms. 
For example, in Equation 
\ref{hllbalx} for multidimensional flux $\mathbf{F_2}^*_{i+\frac12,j+\frac12}$, the third term on right hand side consist of terms $\mathbf{{G_2}}_{RU}$, $\mathbf{{G_2}}_{LU}$, 
$\mathbf{{G_2}}_{LD}$ and $\mathbf{{G_2}}_{RD}$ clearly indicating the influence of y-directional flux terms in the evaluation of x-directional
flux. Equation \ref{hllbaly} depicts a vice versa situation for the flux $\mathbf{G_2}^*_{i+\frac12,j+\frac12}$.
These multidimensional coupling terms along with the wave averaged convective terms may be providing the additional cross dissipation that damps any unprecedented growth in 
instabilities thereby making the present scheme immune to shock instabilities. These equations
reveal that the discretization of convective and pressure fluxes in the corner of
interfaces using GM-HLLE strategy may have a positive impact in curing such instabilities.
\section{Conclusions}
\label{conclusion}

The present work introduces a new genuinely multidimensional contact preserving Riemann solver GM-K-CUSP-X. This scheme is based upon 
Toro-Vazquez type PDE level flux splitting and the ensuing convective-pressure fluxes are discretized following \cite{mandal_2012}. 
While the convective fluxes are upwinded based on appropriate wave speeds that emerge from the interacting states, the pressure fluxes are treated
in a HLLE framework. Restoration of stationary contact preservation ability is improved by explicitly reducing the numerical dissipation in the pressure flux
discretization. 
The resulting solver is found to produce improved results as compared to the original K-CUSP-X solver on standard test problems. 
Particularly, interesting flow features like Kelvin-Helmholtz roll up and mushroom cap structure in two dimensional Riemann problem and 
complex shock interactions in double Mach reflection problem are well resolved. 
Further, the present genuinely multidimensional solver is able to mitigate various forms of shock instabilities that plagued the corresponding
conventional two state Riemann solver, K-CUSP-X. 
%
Such a finding reassures the pre existing notion that multidimensional dissipation is one of the most promising methods 
to cure shock instabilities.
Due to the simplicity of the formulation, the present solver can be easily extended to unstructured framework and three dimensional problems in principle.

\bibliographystyle{ieeetr}

\end{document}